\documentclass[conference]{IEEEtran}
\usepackage{cite}
\usepackage{amsmath,amssymb,amsfonts}
\usepackage{algorithmic}
\usepackage{graphicx}
\usepackage{textcomp}
\usepackage{xcolor}
\usepackage{graphicx}

\def\BibTeX{{\rm B\kern-.05em{\sc i\kern-.025em b}\kern-.08em
    T\kern-.1667em\lower.7ex\hbox{E}\kern-.125emX}}

\newtheorem{definition}{Definition}[section]
\begin{document}

\title{Visually Representing the Landscape of Mathematical Structures}

\author{\IEEEauthorblockN{Katherine Gravel}
\IEEEauthorblockA{\textit{Department of Mathematics} \\
\textit{Massachusetts Institute of Technology}\\
Cambridge, Massachusetts \\
kgravel@mit.edu}
\and
\IEEEauthorblockN{Hayden Jananthan}
\IEEEauthorblockA{\textit{Department of Mathematics}\\
\textit{Vanderbilt University}\\
Nashville, Tennessee\\
hayden.r.jananthan@vanderbilt.edu}
\and
\IEEEauthorblockN{Jeremy Kepner}
\IEEEauthorblockA{\textit{Lincoln Laboratory Supercomputing Center} \\
\textit{Massachusetts Institute of Technology}\\
Lexington, Massachusetts \\
kepner@ll.mit.edu}
}

\maketitle

\begin{abstract}
The information technology explosion has dramatically increased the application of new mathematical ideas and has led to an increasing use of mathematics across a wide range of fields that have been traditionally labeled ``pure'' or ``theoretical''.  There is a critical need for tools to make these concepts readily accessible to a broader community. This paper details the creation of visual representations of mathematical structures commonly found in pure mathematics to enable both students and professionals to rapidly understand the relationships among and between mathematical structures. Ten broad mathematical structures were identified and used to make eleven maps, with 187 unique structures in total. The paper presents a method and implementation for categorizing mathematical structures and for drawing relationship between mathematical structures that provides insight for a wide audience. The most in depth map is available online for public use \cite{b0}.
\end{abstract}

\renewcommand\IEEEkeywordsname{Keywords}

\section{Introduction}
\let\thefootnote\relax\footnotetext{This material is based in part upon work supported by the NSF under grant number DMS-1312831.  Any opinions, findings, and conclusions or recommendations expressed in this material are those of the authors and do not necessarily reflect the views of the National Science Foundation.}
It is often helpful to package information into visuals, which are easily understood \cite{b1}. Visuals tools are used to enhance understanding in media, education, professional circles, and almost every other field. Many times, one may find it more difficult to understand nonvisual media, especially when media is attempting to convey complex or spatial information, such as a direction manual for setting up a couch that does not come with pictures. Visuals representations are also of great use in education, specifically in mathematics, where they are able to ground abstract concepts in a more tangible medium: the physical world \cite{b2,b3,b4,b5}. 

In mathematics as a whole, there are many visual intuitions and tools, although they are concentrated around less abstract mathematics. In grade school, one is taught to count on one's fingers, to ground the abstract concept of quantity in something physical. In middle school and high school, graphs are relied upon heavily to aid in the understanding of slope as rise over run or representing integrals as the area under a curve. Even in undergraduate years, visual intuition still arises. For example, in linear algebra, matrix multiplication is represented by linear transformations. However, as mathematics students move from basic undergraduate courses to higher level subjects including analysis, abstract algebra, and topology, visual intuition and tools for understanding still exist but become much sparser. 

The lack of abundant, clear, visual explanations in abstract mathematics occurs naturally since abstract mathematics seeks to remove any dependency on the physical world and strives to make logical statements based only on generally stated rules \cite{b6}. Since the human brain is wired to process information visually, it is difficult to see how abstract math and visual intuition may intersect \cite{b7,b8}. Therefore, there is a need to create a representative visual that is both abstracted, to mirror the nature of pure mathematics, and specific, so individuals may be able to associate objects represented in the map with something tangible \cite{b9,b10}. 

Visual tools with the above mentioned ideas in mind have been created prior. For example, some individuals, generally students, make small concept maps of different structures in specific mathematical disciplines, or make concept maps of a group of related disciplines and exchange them on internet forums \cite{b11,b12}. Although those maps may be useful to the individuals who created them, they are of less utility to learning communities since they generally do not contain more than four words per node and do not explain the relationships between included concepts. There are also more rigorous representations such as a project by Wolfram Alpha, that allows users to enter a mathematical structure into a search function that will return properties of the object, as well as a local map of what it is related to \cite{b13}. Furthermore, in the book \textit{Mathematics of Big Data} there is extensive use of depictions of the abstract mathematics necessary to define an associative array from more basic mathematical structures \cite{b14}. 

Drawing inspiration from the graphics in \textit{Mathematics of Big Data}, this project created a similar representations of abstract structures while increasing the scope to a large group of fundamental structures in abstract mathematics. This project produced a catalog of different maps representing the connections between abstract mathematical objects. It also categorizes objects based on properties, gives a list of axioms sufficient to constructing each structure and includes other features that make the map aesthetically appealing.  

\section{Background and Definitions} 
Informally, a \emph{mathematical obejct} is any structure built from existing mathematical objects and structures, beginning with sets as the basic building block. Matrices, functions, relations, and sets are all examples of mathematical objects. The terms ``object''€™ and ``€˜structure''€™ will be used interchangeably.  For clarity, we define what is meant by certain terms central to the project:

\begin{definition}{Relation}
An \emph{n-ary relation} $R$ between sets $X_1, \dots,X_n$ is a subset of $X_1 \times \cdots \times X_n$
\label{def: relation}
\end{definition}

\begin{definition} {Function}
An \emph{$n$-ary function} $f: X_0 \times \cdots \times X_n \to Y$ is an $(n+1)$-ary relation $f \subseteq X_0 \times \cdots \times X_n \times Y$ where for each $x_0 \in X_0,\ldots, x_n \in X_n$ there exists a unique $y \in Y$ such that $(x_0,\ldots,x_n,y) \in f$. Such unique $y$ is denoted by $f(x_0,\ldots, x_n)$
\label{def: function}
\end{definition}

\section{Approach}
The first step was to survey pure mathematics as a whole, selecting objects that were central enough to the entire field that most pure mathematicians would know them, excluding categories, which were too board to fall within the scope of the project. The research involved a broad range of textbooks and internet resources \cite{b14,b15,b16,b17,b18}. The next step was to determine a hierarchical relationship between all of the structures. The hierarchical relationships were created by determining the amount of structure encoded by the axioms of an object. For example, there is more structure encoded in the axioms of a module than those of a group. 

The overall layout of the objects were facilitated with the concept map program CMAP \cite{b19}.  CMAP allows for rapid  manipulation of visuals.  Nodes were written to represent each object and labeled arrows were drawn between the nodes to represent their relationship to each other. Although all mathematical objects are related to all other mathematical objects, drawing most connections would be redundant and clutter the diagram. Therefore, the connections that were drawin on maps were unique (with some exceptions disccused later in the paper) Much of the thought of how the connections should be drawn and what the placement of each structure should encode was done at this stage.

The objects were next enumerated in LaTeX as boxes and given descriptions that were rigorous and as complete as possible while preserving a high degree of simplicity. Information about each structure was entered in its corresponding node in a uniform and consistent manner. Using TikZ \cite{b20}, each node was positioned using the placement determined in the CMAP stage. Also, with the assistance of TikZ, arrows were drawn between relevant structures and labeled with brief descriptions of the relationship between the structures. The overall map was also shaded so that it would be easier for an individual to use. 

\section{Design}
Balancing concise and complete information in the diagrams was very important for user readability. 
There were four main categories that all characteristics of an object could be described as: 
\begin{description}
\item[Types] \hfill \\ Previously defined objects.
\item[Functions]\hfill \\ Described functions (Definition~\ref{def: function}) that must be defined over the structure.
\item[Relations]\hfill \\ Enumerated the relations ((Definition~\ref{def: relation}) that must be defined over the structure.
\item[Properties]\hfill \\ Detailed additional properties that functions and relations over the structure must satisfy. 
\end{description} 

For ease of presentation, the functions and relations on the structures were not written in terms of set theoretic functions (e.g. union, intersection) and relations (e.g. membership). Although set theoretic definitions could be written in principle to ensure that the functions and relations over structures really behaved as they claimed to, it would be cumbersome and hinder the definitions that are used more commonly in literature and most fields of mathematics.  

Many structures included in the diagrams have multiple accepted and useful definitions. In most cases, it did not make sense to list more than one definition, since there was not a consistent or elegant way to do so. In general, the definition that would be presented earliest to a student learning in a traditional classroom. However, often the relationship between structures it was related to would be defined normally in terms of an equivalent definition, so for consistency, the relationship would need to be modified to be in terms of the definition stated within the node. The only exception depicting on the maps were Lattices, which have an algebraic and an order theory defintion. Niether could be excluded because many of the different types of lattices, both order theoretic and algebraic properties of the lattice were used in definitions.  

While some structures have multiple accepted standard definitions, there are other structures that have no standard definition. For example, some but not all authors require that rings have a multiplicative identity. To deal with nonstandard definitions, as few choices for definitions were made as possible. Although it would be possible to include multiple definitions, it was decided that it would quickly become cumbersome. When choices in definitions were made, they were made to follow the what the majority of authors defined the structure as. Distinctions generally only needed to be made at a low level of complexity, as structures that built upon them only needed to reference the structure as a whole, and not the internal mechanisms of the structure. Additionally, the nuances of operations were not stated in the diagrams, which is where most of subtleties involving non standard definitions applied. For example, a module is an abelian group that is closed under ring multiplication. But the mechanism of scalar multiplication is not discussed, which is where it is determined what properties a ring of scalars must satisfy, so there is no conflict. 

Precedence was given to simplicity and intuitive representation. Therefore, certain technically unnecessary characteristics of structures were only listed if they gave a significantly fuller picture of the structure. For example, when describing a lie algebra, it is redundant to note that both $[X,X] = 0$ and $[X,Y] = -[Y,X]$, but was included because both give insight to the nature of a lie bracket. 

Although there are many fundamental relationships between different mathematical structures (e.g. the association of a Lie group with its corresponding Lie algebra), for simplicity and to avoid determining what constitutes being `fundamental' in general, only associations which indicate one structure \emph{extending} another by the additional of new functions, relations, or properties were included. The inclusion of examples of structures was on a case-by-case basis, due in large part to the subtlety of distinguishing between a class of examples of a particular structure and a new mathematical structure. The mathematical literature as well as relationship with other mathematical structure served as a guide for this distinction. For example, General Linear Groups are not included while Cyclic Groups are. 

The hierarchical structure of each map was further enhanced as follows.
\begin{description}
\item[Size of Nodes] \hfill \\ Nodes that were `higher' on the hierarchy, meaning they were more general and less specific (e.g. set), were given larger area while more specific and structured nodes (e.g. primitive group) were given smaller dimensions. 
\item[Coloring] \hfill \\ Different sections of structures were colored differently. Ten large sections were identified (algebras, fields, graphs, groups, lattices, posets, modules, rings, sets, topological spaces). Each was given a distinct color.
\end{description}

Some structures, such as a topological group, should be classified in more than one section: groups and topological spaces. To deal with the coloring issue, nodes were simply shaded multiple colors through fading. In the example of topological group, the topological group node faded from orange (the color for groups) to purple (the color for topological spaces). Finally, all nodes on all maps were hyperlinked to their wikipedia pages, so a user may further explore mathematical concepts disscussed relatively informally.  

\section{Results}
\begin{figure}
	\includegraphics[width=\linewidth]{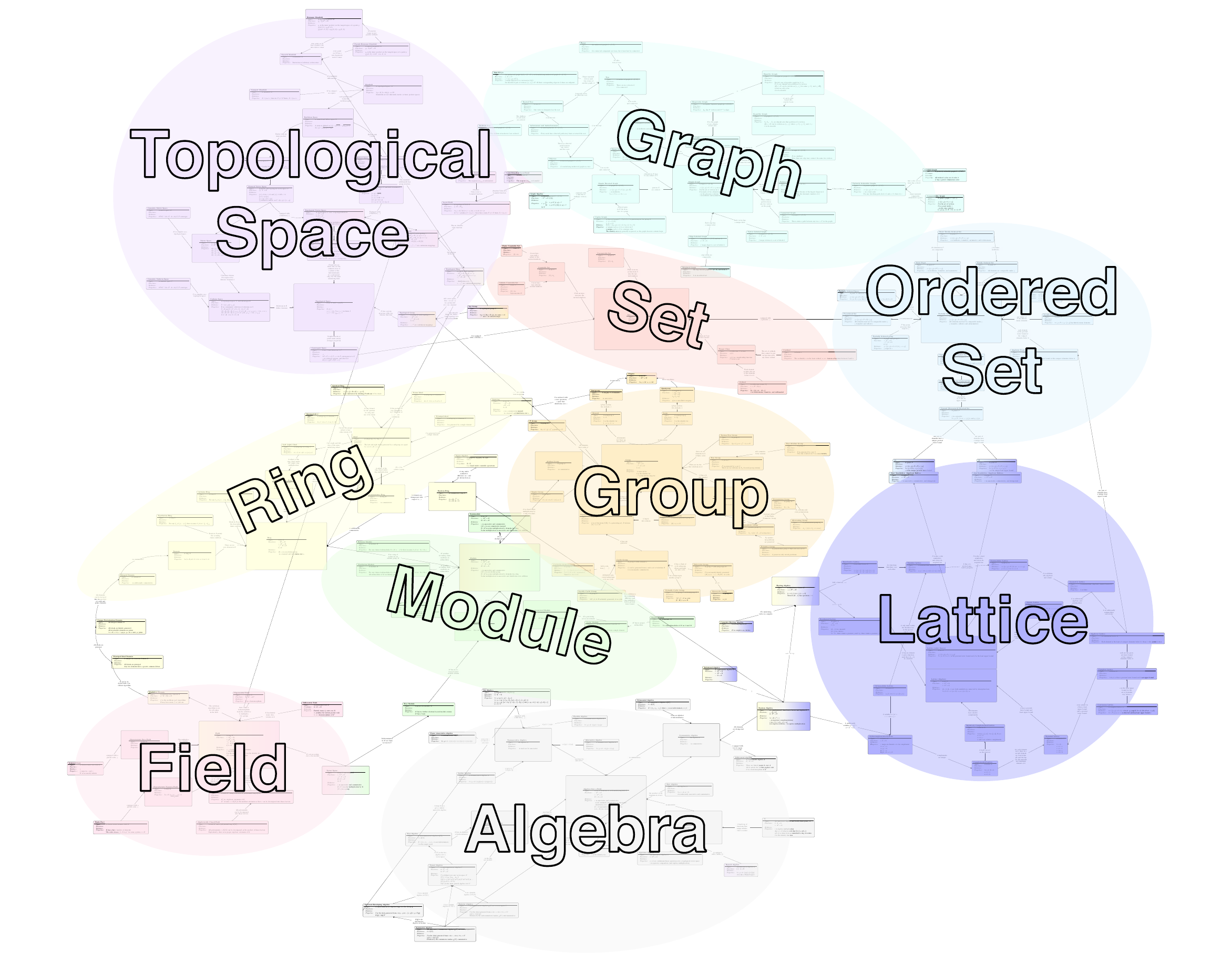}
	\caption{Map of the largest mathematical structures.}
	\label{fig:BigMap}
\end{figure}

This project created eleven maps, totaling 187 unique structures. The map shown in Figure~\ref{fig:BigMap} contained only all of the structures. The overlay labels the large catagories of structures that were included. Nine maps contained all of the related structures grouped under single color (with the exception of posets and lattices, which were combined into one map with two colors). The final map contained only the largest structures, meaning the structures that recieved their own color. 

\begin{figure} [h!]
	\includegraphics[width=\linewidth]{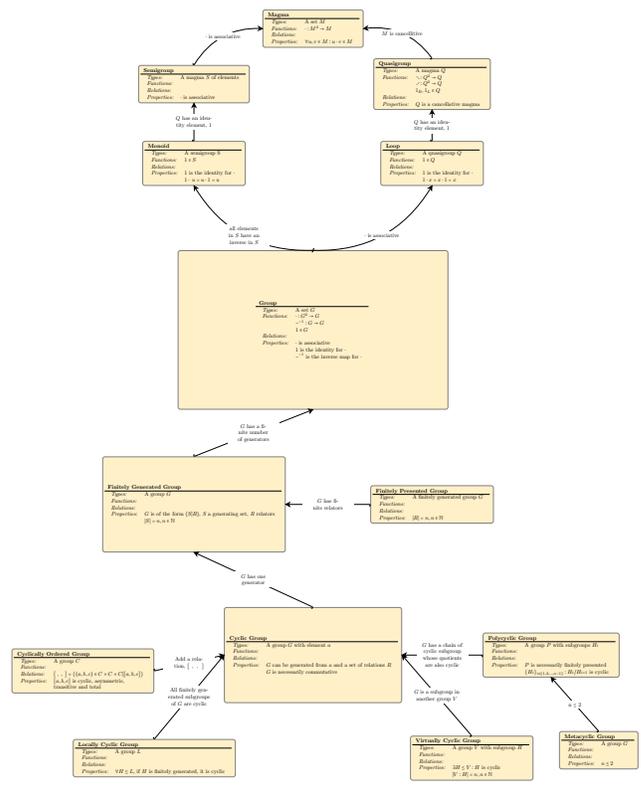}
	\caption{An outline of connected structures.}
	\label{fig:outline}
\end{figure}

\ref{fig:outline} depicts the structure of a simplified map of groups. It is noted that all nodes have information about the structures they represent, nodes are connected as minimally as possible, nodes are sized differently based on hierarchy, and arranged strategically to encode information and all nodes are purposefully colored orange, the color designated to groups.

\begin{figure} [h!]
	\includegraphics[width=\linewidth]{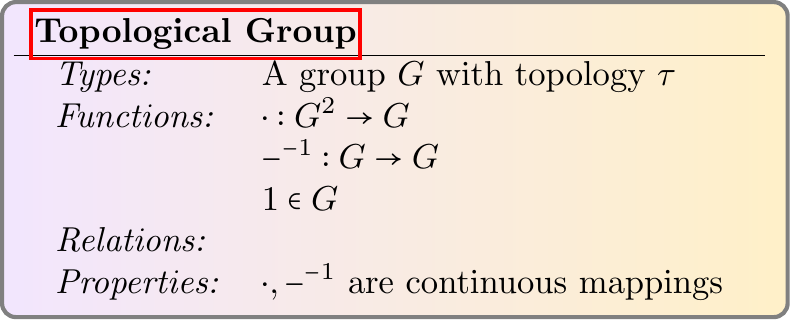}
	\caption{An individual node. It has types of group and topology, so t is shaded both orange, in accordance with the group color scheme, and purple, to match topological space coloring. }
	\label{fig:singalnode}
\end{figure}

This project was designed to simplify and accelerate the understanding of how abstract structures, stripped of their commonly placed disciplines, interact with and relate to each other. The use of labeled arrows between structures is meant to clearly depict the relationship they had with each other, and the use of consistent definitions is meant to make the distinctions between the structures as clear as possible. Aesthetic features, such as coloring, were used to highlight important subtleties that arise from certain structures. Furthermore, a hierarchical structure was also created in these maps, which gives individuals a framework to organize definitions of related but separated structures in their mind through the use of differently sized nodes and inclusion of structures in the \textit{Types} of other structures. 

\section{Summary}
The project created a unique representation of common mathematical objects aimed at young mathematical students and professionals in related fields with none to minimal background in pure mathematics. Eleven maps were created in total, with careful attention paid to the visual details of the map, as conveying a visual representation of abstract concepts was especially important in a field that lacks ample visual representations. A single map with all 187 nodes is posted online and available for use \cite{b0}.

\section*{Acknowledgments}
The authors wish to acknowledge the following individuals for their contributions and support:
Alan Edelman, Vijay Gadepally, Chris Hill, and Lauren Milechin.

\end{document}